\documentclass[12pt]{article}
\usepackage{e-jcb}
\usepackage{amsmath}
\usepackage{amssymb}
\usepackage{amsthm}
\usepackage{graphicx}
\usepackage[all]{xy}
\usepackage[square, numbers, longnamesfirst]{natbib}
\usepackage{caption}
\usepackage{subcaption}
\usepackage[colorlinks=true,citecolor=black,linkcolor=black,urlcolor=blue]{hyperref}

\newcommand{\doi}[1]{\href{http://dx.doi.org/#1}{\texttt{doi:#1}}}

\DeclareMathOperator{\cone}{cone}

\def\vec#1{\mathchoice{\mbox{\boldmath$\displaystyle\bf#1$}}
{\mbox{\boldmath$\textstyle\bf#1$}}
{\mbox{\boldmath$\scriptstyle\bf#1$}}
{\mbox{\boldmath$\scriptscriptstyle\bf#1$}}}
\DeclareMathOperator{\sign}{sgn}

\theoremstyle{plain}
\newtheorem{theorem}{Theorem}

\theoremstyle{definition}
\newtheorem{definition}[theorem]{Definition}
\newtheorem{example}[theorem]{Example}
\newtheorem{conjecture}[theorem]{Conjecture}

\theoremstyle{remark}
\newtheorem{remark}[theorem]{Remark}

\numberwithin{theorem}{section}

\title{\bf The unreasonable ubiquitousness of quasi-polynomials\thanks{With apologies to \citet{W60} and \citet{H80}. Extended abstract appeared in FPSAC 2013. This version to appears in the Electronic Journal of Combinatorics.}}
\author{Kevin Woods\\
\small Department of Mathematics\\[-0.8ex]
\small Oberlin College\\[-0.8ex] 
\small Oberlin, Ohio, USA\\
\small\tt Kevin.Woods@oberlin.edu\\
}

\date{\small Mathematics Subject Classifications: 05A15, 52C07, 03B10}

\begin{document}

\providecommand{\abs}[1]{\lvert#1\rvert}
\providecommand{\floor}[1]{\left\lfloor#1\right\rfloor}
\providecommand{\Z}{\mathbb{Z}} \providecommand{\R}{\mathbb{R}}
\providecommand{\N}{\mathbb{N}} \providecommand{\C}{\mathbb{C}}
\providecommand{\Q}{{\mathbb{Q}}} \providecommand{\x}{\mathbf{x}}
\providecommand{\y}{\mathbf{y}} \providecommand{\z}{\mathbf{z}}

\maketitle
\begin{abstract}
  A function $g$, with domain the natural numbers, is a quasi-poly\-no\-mial if there exists a period $m$ and polynomials $p_0,p_1,\ldots,p_{m-1}$ such that $g(t)=p_i(t)$ for $t\equiv i\bmod m$. Quasi-polynomials classically -- and ``reasonably'' -- appear in Ehrhart theory and in other contexts where one examines a family of polyhedra, parametrized by a variable $t$, and defined by linear inequalities of the form $a_1x_1+\cdots+a_dx_d\le b(t)$.
  
  Recent results of Chen, Li, Sam; Calegari, Walker; and Roune, Woods show a quasi-polynomial structure in several problems where the $a_i$ are also allowed to vary with $t$. We discuss these ``unreasonable'' results and conjecture a general class of sets that exhibit various (eventual) quasi-polynomial behaviors: sets $S_t\subseteq\N^d$ that are defined with quantifiers ($\forall$, $\exists$), boolean operations (and, or, not), and statements of the form $a_1(t)x_1+\cdots+a_d(t)x_d \le b(t)$, where $a_i(t)$ and $b(t)$ are polynomials in $t$. These sets are a generalization of sets defined in the Presburger arithmetic. We prove several relationships between our conjectures, and we prove several special cases of the conjectures. The title is a play on Eugene Wigner's ``The unreasonable effectiveness of mathematics in the natural sciences''.

 \bigskip\noindent \textbf{Keywords:} Ehrhart polynomials, generating functions, Presburger arithmetic, quasi-polynomials, rational generating functions
\end{abstract}

\section{Reasonable Ubiquitousness}
\label{sec1}
In this section, we survey classical appearances of quasi-polynomials (though Section \ref{subsec:PA} might be new even to readers already familiar with Ehrhart theory). In Section \ref{sec2}, we survey some recent results where the appearance of quasi-polynomials is more surprising. In Section \ref{sec3}, we make several conjectures generalizing these ``unreasonable'' results. We state theorems relating these conjectures and state theorems proving certain cases. In particular, we conjecture that any family of sets $S_t$ -- defined with quantifiers ($\forall$, $\exists$), boolean operations (and, or, not), and statements of the form $\vec a(t)\cdot \x \le b(t)$ (where $\vec a(t)\in\Z[t]^d, b(t)\in \Z[t]$, and $\cdot$ is the standard dot product) -- exhibits eventual quasi-polynomial behavior, as well as rational generating function behavior. 
Of course, reasonable people may disagree on what is unreasonable; the title is a play on Eugene Wigner's ``The unreasonable effectiveness of mathematics in the natural sciences'' \cite{W60}. All proofs are contained in Section \ref{sec4}. We use bold letters such as $\vec x=(x_1,\ldots,x_d)$ to indicate multi-dimensional vectors.

\begin{definition}
A function $g:\N\rightarrow \Q$ is a \emph{quasi-polynomial} if there exists a period $m$ and polynomials $p_0,p_1,\ldots,p_{m-1}\in\Q[t]$ such that
\[g(t)=p_i(t),\text{ for }t\equiv i\bmod m.\]
\end{definition}

\begin{example}
\label{ex:floor}
\[g(t)=\floor{\frac{t+1}{2}}=\begin{cases} \frac{t}{2} & \text{if $t$ even},\\ \frac{t+1}{2} & \text{if $t$ odd},\end{cases}\]
is a quasi-polynomial with period 2.
\end{example}

This example makes it clear that the ubiquitousness of quasi-polynomials shouldn't be  too surprising: anywhere there are floor functions, quasi-poly\-no\-mials are likely to appear.

Floor functions can also create multivariate quasi-polynomials:

\begin{example}
\label{ex:mfloor}
Let $g(s,t)=\floor{(t+s+1)/2}$. The behavior of $g$ depends on the parity of $t+s+1$. Let $\Lambda$ be the lattice $(0,2)\Z+(1,1)\Z$, the set of $(s,t)\in\Z^2$ for which $t+s+1$ is odd. Then the lattice coset $(0,1)+\Lambda$ is the set for which $t+s+1$ is even, and 
\[g(s,t)=\begin{cases}\frac{t+s}{2} &\text{if $(s,t)\in\Lambda$},\\ \frac{t+s+1}{2} &\text{if $(s,t)\in (0,1)+\Lambda$}.\end{cases}\]
\end{example} 

This motivates the definition of a multivariate quasi-polynomial; ``residues modulo a period'' are replaced by ``cosets of a lattice'':

\begin{definition}
$g:\N^n\rightarrow\Q$ is a (multivariate) \emph{quasi-polynomial} if there exists an $n$-dimensional lattice $\Lambda\subseteq\Z^n$, a set $\{\vec \lambda_i\}$ of coset representatives of $\Z^n/\Lambda$, and polynomials $p_i\in\Q[\vec t]$ such that
\[g(\vec t)=p_i(\vec t),\text{ for }\vec t\in \vec \lambda_i+\Lambda.\]
\end{definition}


We will generally be concerned with \emph{integer-valued} quasi-poly\-no\-mials, those quasi-poly\-no\-mials whose range lies in $\Z$. Note that Examples \ref{ex:floor} and \ref{ex:mfloor}
demonstrate that such quasi-polynomials may still require rational coefficients.

\subsection{Ehrhart theory}
Perhaps the most well-studied quasi-polynomials are the \emph{Ehrhart quasi-poly\-no\-mials}:

\begin{theorem}[\citet{Ehrhart62}]
\label{thm:ehr}
Suppose $P$ is a polytope (bounded polyhedron) whose vertices have rational coordinates. Let $g(t)$ be the number of integer points in $tP$, the dilation of $P$ by a factor of $t$. Then $g(t)$ is a quasi-polynomial, with period the smallest $m$ such that $mP$ has integer coordinates.
\end{theorem}

\begin{example}
\label{ex:square}
Let $P$ be the triangle with vertices $(0,0)$, $\left(\frac{1}{2},0\right)$, and $\left(\frac{1}{2},\frac{1}{2}\right)$. Then
\begin{equation}\label{eq:Eh}g(t)=\#(tP\cap\Z^2)=\frac{\big(\floor{t/2}+1\big)\big(\floor{t/2}+2\big)}{2}=\begin{cases} (t+2)(t+4)/8&\text{if $t$ even,}\\ (t+1)(t+3)/8 & \text{if $t$ odd,}\end{cases}\end{equation}
is a quasi-polynomial with period 2.
\end{example}

Writing $tP$ from this example as
\[\big\{(x,y)\in\R^2:\ 2x\le t,\ y-x\le 0,\ -y\le 0\big\}\]
suggests a way to generalize this result: for $\vec t\in\N^n$, let $S_{\vec t}$ be the set of integer points, $\x\in\Z^d$, in a polyhedron defined with linear inequalities of the form $\vec a \cdot \vec x\le b(\vec t)$, where $\vec a\in\Z^d$ and $b(\vec t)\in \Z[\vec t]$ is linear.

\begin{figure}
\begin{center}
\begin{tabular}{c|c|c}
\includegraphics[width=1.5in]{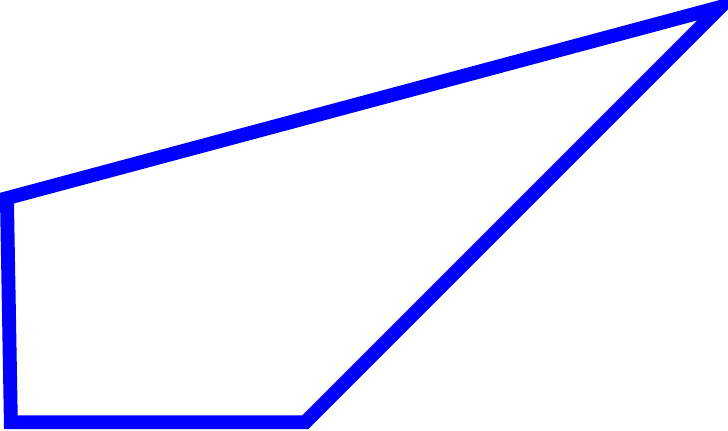} &
\includegraphics[width=1.5in]{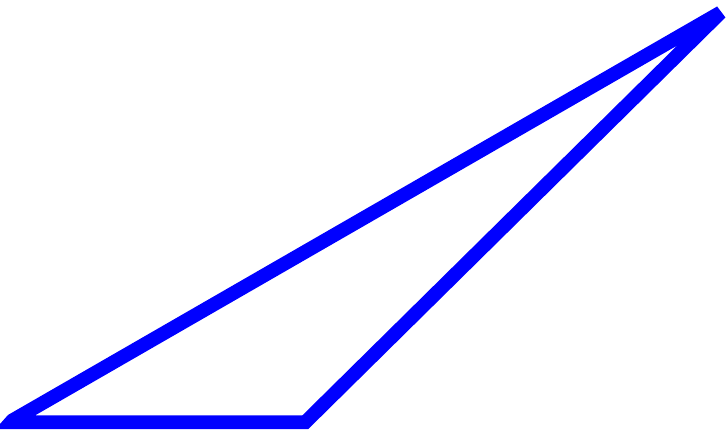} &
\includegraphics[width=1.5in]{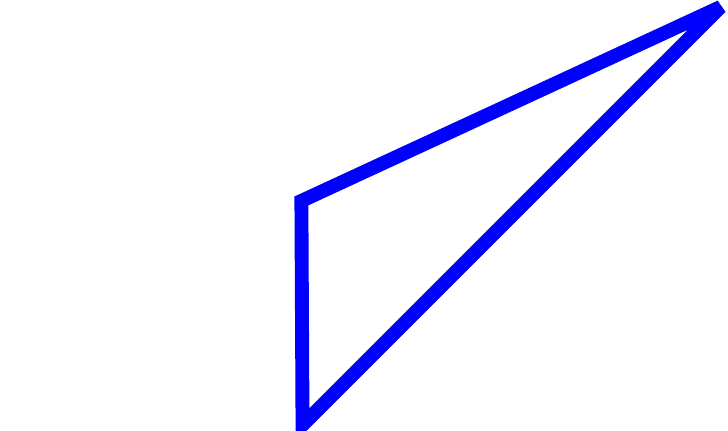} \\
$t\le s\le2t$ & $0\le 2t\le s$ & $0\le s\le t$
\end{tabular}
\caption{Polyhedra defined in Example \ref{ex:sven} for various $(s,t)\in\N^2$.}
\label{fig:sven}
\end{center}
\end{figure}

\begin{example}
\label{ex:sven}
Let 
\[S_{s,t}=\big\{(x,y)\in\Z^2:\ 2y-x\le 2t-s,\ x-y\le s-t,\ x,y\ge 0\big\}.\]
\end{example}
For a fixed $(s,t)$, $S_{s,t}$ is the set of integer points in a polyhedron in $\R^2$. As $(s,t)$ varies, the ``constant'' term of these inequalities change, but the coefficients of $x$ and $y$ do not; in other words, the normal vectors to the facets of the polyhedron do not change, but the facets move ``in and out''. In fact, they can move in and out so much that the combinatorial structure of the polyhedron changes. Figure \ref{fig:sven} shows the combinatorial structure for different $(s,t)\in\N^2$. Using various methods, \citet{Beck04} and \citet{VW04} compute that

\[g(s,t)=\abs{S_{s,t}}=\begin{cases}
\frac{s^2}2 - \lfloor \frac{s}2 \rfloor s + \frac {s} 2 + \lfloor \frac{s}2 \rfloor^2 + \lfloor \frac{s}2 \rfloor +1 & \text{if }t\le s\le2t,\\
s t - \lfloor \frac{s}2 \rfloor s - \frac {t^2}2 + \frac {t} 2+ \lfloor \frac{s}2 \rfloor^2 + \lfloor \frac{s}2 \rfloor +1 & \text{if }0\le 2t\le s,\\
\frac {t^2}2 + \frac {3t}2 +1 & \text{if }0\le s\le t.
\end{cases}
\]

In this example, the function $g(s,t)$ is a quasi-polynomial, at least on pieces of parameter-space for which the combinatorial type is constant.

\begin{definition}
A function $g:\N^n\rightarrow \Q$ is a \emph{piecewise quasi-polynomial} if there exists a finite
partition $\bigcup_i (P_i\cap\N^n)$ of $\N^n$ with $P_i$ polyhedra (which may not all be full-dimensional)
and there exist quasi-polynomials $g_i$ such that
\[g(\vec t)=g_i(\vec t) \text{ for }\vec t\in P_i\cap\N^n.\]  
\end{definition}

\citet{Sturmfels95} effectively proved the following generalization of Ehrhart theory:

\begin{theorem}
\label{thm:par_polyhedra}
Let $S_{\vec t}$ be the set of integer points, $\x\in\Z^d$, in a polyhedron defined with linear inequalities of the form $\vec a \cdot \vec x\le b(\vec t)$, where $\vec a\in\Z^d$ and $b(\vec t)\in\Z[\vec t]$ is linear. Then $g(\vec t)=\abs{S_{\vec t}}$ is a piecewise quasi-polynomial.
\end{theorem}

\begin{remark}
Sections \ref{sec2} and further will predominantly be concerned with \emph{univariate} functions. Being a univariate piecewise quasi-polynomial $g:\N\rightarrow \Q$ is equivalent to \emph{eventually} being a quasi-polynomial; that is, there exists a $T$ such that for all $t\ge T$, $g(t)$ agrees with a quasi-polynomial.
\end{remark}

\subsection{Generating functions}
\label{subsec:gfs}
Many classical proofs of Ehrhart's Theorem (Theorem \ref{thm:ehr}) use generating functions. To prove that a function $g(t)$ is a quasi-polynomial of period $m$, it suffices (see Section 4.4 of \citet{Stanley12}) to prove that the Hilbert series $\sum_{t\in\N}g(t)y^t$ can be written as a rational function of the form
\[\frac{p(y)}{(1-y^m)^d},\]
where $p(y)$ is a polynomial of degree less than $md$. For $g(t)=\abs{tP\cap\Z^2}$ with $P$ the triangle in Example \ref{ex:square}, we can see that
\begin{equation}\label{eq:HS}\sum_{t\in\N}g(t)y^t=1+y+3y^2+3y^3+6y^4+\cdots=\frac{1+y}{(1-y^2)^3}.\end{equation}
Indeed, these proofs of Ehrhart's Theorem start by considering the generating function $\sum_{t\in\N, \vec s\in tP\cap\Z^d}\x^{\vec s}y^t$ (where $\x^{\vec s}=x_1^{s_1}\cdots x_d^{s_d}$) and substituting in $\x=(1,\ldots,1)$ to get the Hilbert series. For $P$ in Example \ref{ex:square},
\begin{align*}
\sum_{t\in\N,\vec s\in tP\cap\Z^2}\hspace{-4pt}\x^{\vec s}y^t&=1+y+(1+x_1+x_1x_2)y^2+(1+x_1+x_1x_2)y^3+(1+\cdots+x_1^2x_2^2)y^4
+\cdots\\
&=\frac{1+y}{(1-y^2)(1-x_1y^2)(1-x_1x_2y^2)}, 
\end{align*}
as can be checked by expanding as a product of infinite geometric series. Substituting $x_1=x_2=1$ yields the Hilbert series in (\ref{eq:HS}).

\begin{definition}
We call any generating function or Hilbert series a \emph{rational generating function} if it can be written in the form
\[\frac{p(\x)}{(1-\x^{\vec b_1})\cdots(1-\x^{\vec b_k})},\]
where $p$ is a Laurent polynomial over $\Q$ and $\vec b_i\in\Z^d$ are lexicographically positive (first nonzero entry is positive).
\end{definition}

While we will generally be assuming that the generating functions are for subsets of $\N^d$, we need $\vec b_i$ to be lexicographically positive rather than simply in $\N^d\setminus\{0\}$ for examples like the following:

\begin{example}
\label{ex:negs}
Let $S=\big\{(x,y)\in\N^2:\ x+y=1000\big\}$. While $y^{1000}+xy^{999}+\cdots+x^{1000}$ is a legitimate rational generating function, it makes more sense to write it as
\[\frac{y^{1000}-x^{1001}y^{-1}}{1-xy^{-1}}.\]
\end{example}

\begin{remark}\label{rmk:lexpos2}
Having $\vec b$ lexicographically positive guarantees that $1/(1-\vec x^{\vec b})=1+\vec x^{\vec b}+\vec x^{2\vec b}+\cdots$ is the Laurent series convergent on a neighborhood of $\vec a \doteq (e^{-\varepsilon},e^{-\varepsilon^2},\ldots,e^{-\varepsilon^d})$, for sufficiently small $\varepsilon$, as follows: Let $b_i>0$ be the first nonzero coordinate of $\vec b$. Then
\[\ln\left(\vec a ^ {\vec b}\right)=-\varepsilon b_1-\varepsilon^2 b_2-\cdots -\varepsilon^d b_d=-\varepsilon^i b_i-\cdots-\varepsilon^d b_d.\]
For sufficiently small $\varepsilon$, the  $-\varepsilon^i b_i$ term dominates, and $\ln(\vec a ^ {\vec b})<0$. This implies that $\abs{\vec a ^ {\vec b}}<1$, and indeed $1+\vec x^{\vec b}+\vec x^{2\vec b}+\cdots$ does converge on a neighborhood of $\vec a$.
\end{remark}

\begin{remark}\label{rmk:lexpos}
If $\vec b$ is lexicographically negative, then we may instead write
\[\frac{1}{1-\vec x^{\vec b}}=\frac{1}{1-\vec x^{\vec b}}\cdot\frac{-\vec x^{-\vec b}}{-\vec x^{-\vec b}}=\frac{-\vec x^{-\vec b}}{1-\vec x^{-\vec b}}\]
with $-\vec b$ lexicographically positive.
\end{remark}

Several classic results start with a simple generating function, and use it to find quasi-polynomial behavior.
\begin{definition}
\label{def:vpf}
Given $\vec a_{1},\ldots,\vec a_{d}\in\N^{n}$, the \emph{vector partition function} $g:\N^{n}\rightarrow\N$ is defined by
\[g(\vec t)=\#\{(\lambda_{1},\ldots,\lambda_{d})\in\N^{d}:\ \vec t=\lambda_{1}\vec a_{1}+\cdots+\lambda_{d}\vec a_{d}\},\]
that is, the number of ways to partition the vector $\vec t$ into parts taken from $\{\vec a_i\}$.
\end{definition}

The generating function $\sum_{\vec t\in\N^n} g(\vec t)\y^{\vec t}$ can be written as
\[\sum_{\vec t\in\N^n} g(\vec t)\y^{\vec t} = \frac{1}{(1-\y^{\vec a_1})\cdots(1-\y^{\vec a_d})},\]
obtained by rewriting the rational function as a product of infinite geometric series. The following is proved by \citet{Sturmfels95}:

\begin{theorem}
\label{thm:vpf}
Any vector partition function is a piecewise quasi-polynomial.
\end{theorem}

See \citet{Beck04} for a self-contained explanation utilizing the partial fraction expansion of the rational function. 
For example, if $a_1=1$ and $a_2=a_3=2$, then the vector partition function is encoded by the generating function
\[\frac{1}{(1-y)(1-y^2)^2}=\frac{1+y}{(1-y^2)^3}.\]
We saw previously that this generating function corresponds to the quasi-polynomial in Example \ref{ex:square}.

In Section \ref{sec3}, we will use a different generating function: for \emph{fixed} $t$, examine the generating function $\sum_{\vec s\in tP\cap\Z^d} \x^{\vec s}$.

\begin{example}In the triangle from Example \ref{ex:square}, this gives us
\[\left(1+x_1+x_1^2+\cdots+x_1^{\floor{t/2}}\right)+\left(x_1+x_1^2+\cdots+x_1^{\floor{t/2}}\right)x_2+\cdots +\left(x_1^{\floor{t/2}}\right)x_2^{\floor{t/2}}.\]
We can write this more compactly as
\begin{equation}\label{eq:gf}\sum_{\vec s\in tP\cap\Z^d} \x^{\vec s}=\frac{1}{(1-x_1)(1-x_1x_2)}-\frac{x_1^{\floor{t/2}+1}}{(1-x_1)(1-x_2)}+\frac{x_1^{\floor{t/2}+1}x_2^{\floor{t/2}+2}}{(1-x_2)(1-x_1x_2)},\end{equation}
which we can verify directly by expanding the fractions as products of geometric series. Given this generating function, we can count the number of integer points in $tP$ by substituting in $\x=(1,\ldots,1)$. Substituting $x_1=x_2=1$ into (\ref{eq:gf}), we see that $(1,1)$ is a pole of these fractions. Fortunately, getting a common denominator and applying L'H\^opital's rule to find the limit as $x_1$ and $x_2$ approach 1 will work, and it is evident that the differentiation involved in L'H\^opital's rule will yield a quasi-polynomial in $t$ as the result; careful calculation will show that it matches (\ref{eq:Eh}).
\end{example}

In general, Proposition 2.11 of \citet{VW04} or Theorem 4.4 of \citet{BP99} both give an algorithmic version of the following theorem:

\begin{theorem}
\label{thm:VW}
For $\vec t\in\Z^n$, let $S_{\vec t}$ be the set of integer points, $\x\in\Z^d$, in a polyhedron defined with linear inequalities of the form $\vec a \cdot \vec x\le b(\vec t)$, where $\vec a\in\Z^d$ and $b(\vec t)\in\Z[\vec t]$ is linear. Then there is a finite decomposition of $\Z^n$ into pieces of the form $P\cap\Z^n$ (with $P$ a polyhedron) such that, considering the $\vec t$ in each piece separately,
\[\sum_{\vec{s}\in S_{\vec t}}\vec{x}^{\vec{s}}=\sum_{i=1}^m\epsilon_i\frac{\vec{ x}^{\vec u_i(\vec t)}}{(1-\vec{x}^{\vec d_{i1}})\cdots(1-\vec{x}^{\vec d_{ik_i}})},\]
where $\epsilon_i=\pm1$, the coordinates of $\vec u_i$ are linear quasi-polynomials in $\vec t$, and $\vec d_{ij}\in\Z^d\setminus\{0\}$.
\end{theorem}

Substituting $\vec x =(1,\ldots,1)$, using L'H\^opital's rule as necessary, recovers Theorem \ref{thm:par_polyhedra}.

A key step in proving Theorem \ref{thm:VW} is Brion's Theorem \cite{Brion88} (see Chapter 9 of \citet{BR07} for a proof and discussion), which reduces computing generating functions for polyhedra to computing generating functions for cones:
\begin{theorem}[Brion's Theorem]
\label{thm:Brion}
Let $P$ be a polyhedron in the nonnegative orthant $\R_{\ge 0}^d$. For each vertex $\vec v$ of $P$, define its \emph{vertex cone} $C_{\vec v}$ to be the set of points in $\R^d$ satisfied by those inequalities defining $P$ that are \emph{equalities} at $\vec v$. Then 
\[\sum_{\vec{s}\in P\cap\Z^d}\vec{x}^{\vec{s}}=\sum_{\vec v\text{ a vertex of $P$}}\ \sum_{\vec s\in C_{\vec v}\cap\Z_d}\vec{x}^{\vec{s}}.\]
\end{theorem}
Note that, for different vertices $\vec v$, the $\sum_{\vec s\in C_{\vec v}\cap\Z_d}\vec{x}^{\vec{s}}$ may converge on different neighborhoods. If they are written as rational functions, Remarks  \ref{rmk:lexpos2} and \ref{rmk:lexpos} fix this problem.

\begin{example}
For $S=\big\{(x,y)\in\N^2:\ x+y=1000\big\}$ from Example \ref{ex:negs}, the two vertex cones are
\[C_1=\big\{(0,1000)+\lambda(1,-1):\ \lambda\ge 0\big\}\quad\text{and}\quad C_2=\big\{(1000,0)+\lambda(-1,1):\ \lambda\ge 0\big\},\]
and so
\begin{align*}\sum_{\vec s\in S}\vec x^{\vec s}=\sum_{\vec s\in C_1\cap\Z^2}\vec x^{\vec s}+\sum_{\vec s\in C_2\cap\Z^2}\vec x^{\vec s}
&=\left(y^{1000}+xy^{999}+\cdots\right)+\left(x^{1000}+x^{999}y+\cdots\right)\\
&=\frac{y^{1000}}{1-xy^{-1}}+\frac{x^{1000}}{1-x^{-1}y}=\frac{y^{1000}-x^{1001}y^{-1}}{1-xy^{-1}},
\end{align*}
as seen in Example \ref{ex:negs}.
\end{example}

\begin{example}
The triangle from Example \ref{ex:square} has three vertex cones. Brion's Theorem applied to $tP\cap\Z^2$ is illustrated in (\ref{eq:gf}), with Remark \ref{rmk:lexpos} already applied.
\end{example}

\begin{remark}\label{rmk:lines}
Note that Brion's Theorem is often stated for polyhedra in $\R^d$ (unrestricted to the nonnegative orthant). If the polyhedron contains straight lines, one must be more careful when talking about generating functions. For example, the set $\Z$ has generating function
\[\cdots+x^{-1}+1+x^1+x^2+\cdots=\frac{x^{-1}}{1-x^{-1}}+\frac{1}{1-x}=-\frac{1}{1-x}+\frac{1}{1-x}=0.\]
See, for example, \citet{Barvinok08} for more details. In this paper, we will constrain the sets that we are interested in to be nonnegative, to avoid issues like this.
\end{remark}

\subsection{Presburger arithmetic}
\label{subsec:PA}

So far, our examples have been integer points in polyhedra. A key property of such sets is that they can be defined without quantifiers. However, even for sets defined with quantifiers, we end up with reasonable appearances of quasi-polynomials.

\begin{definition}
A \emph{Presburger formula} is a boolean formula with variables in $\N$ that can be written using quantifiers ($\exists, \forall$), boolean operations (and, or, not), and linear (in)equalities in the variables. We write a Presburger formula as $F(\vec u)$ to indicate the free variables $\vec u$ (those not associated with a quantifier).
\end{definition}

\citet{Presburger29} (see \citep{Presburger91} for a translation) examined this first order theory and proved it is decidable.

\begin{example}
\label{ex:PA}
Given $t\in\N$, let
\[S_t=\big\{x\in\N:\ \exists y\in\N,\  2x+2y+3=5t\text{ and }t< x\le y\big\}.\]
We can compute that
\[S_t=\begin{cases}\left\{t+1,t+2,\ldots,\floor{\frac{5t-3}{4}}\right\} & \text{if $t$ odd, $t\ge 3$,}\\
\emptyset & \text{else.}
\end{cases}\]
\end{example}

This set has several properties, \textit{cf.} Section \ref{sec3}:
\begin{enumerate}
\item The set of $t$ such that $S_t$ is nonempty is $\{3,5,7,\ldots\}$. This set is eventually periodic.
\item The cardinality of $S_t$ is
\[S_t=\begin{cases}\floor{\frac{5t-3}{4}}-t & \text{if $t$ odd, $t\ge 3$,}\\
0 & \text{else,}
\end{cases}\]
which is eventually a quasi-polynomial of period 4.
\item When $S_t$ is nonempty, we can obtain an element of $S_t$ with the function $x(t)=t+1$, and $x(t)$ is eventually a quasi-polynomial.
\item[3a.] More strongly, when $S_t$ is nonempty, we can obtain the maximum element of $S_t$ with the function $x(t)=\floor{(5t-3)/4}$, and $x(t)$ is eventually a quasi-polynomial.
\item We can compute the generating function
\begin{align*}
\sum_{s\in S_t}x^s&=\begin{cases}x^{t+1}+x^{t+2}+\cdots+x^{\floor{(5t-3)/4}} & \text{if $t$ odd, $t\ge 3$,}\\
0 & \text{else,}
\end{cases}\\
&=\begin{cases}\dfrac{x^{t+1}-x^{\floor{(5t-3)/4)}+1}}{1-x} & \text{if $t$ odd, $t\ge 3$,}\\
0 & \text{else.}
\end{cases}
\end{align*}
We see that, for fixed $t$, this generating function is a rational function. Considering each residue class of $t \bmod 4$ separately, the exponents in the rational function can eventually be written as polynomials in $t$.
\end{enumerate}

Versions of these properties always hold for sets defined in Presburger arithmetic. For example, \citet{Woods13} gave several properties of Presburger formulas that hold even for sets defined with multivariate parameters, $\vec t\in\N^n$:
\begin{theorem}[from Theorems 1 and 2 of \citet{Woods13}]
\label{thm:PA} 
Suppose $F(\vec s, \vec t)$ is a Presburger formula, with $\vec s$ and $\vec t$ collections of free variables. Then
\begin{itemize}
\item $g(\vec t)=\#\big\{\vec s\in\N^d: \ F(\vec s, \vec t)\big\}$ is a piecewise quasi-polynomial,
\item $\sum_{\vec s,\vec t:F(\vec s,\vec t)} \x^{\vec s}\y^{\vec t}$ is a rational generating function, and
\item $\sum_{\vec t\in\N^n}g(\vec t)\y^{\vec t}$ is a rational generating function.
\end{itemize}
\end{theorem}
Property 4 from Example \ref{ex:PA} can be proved in general by using Theorem \ref{thm:PA} to write $\sum_{\vec s,\vec t:F(\vec s,\vec t)} \x^{\vec s}\y^{\vec t}$ as a rational generating function and applying Theorem \ref{thm:diff_gfs}. The proof of Theorem \ref{thm:implications} then shows that all of the other properties follow.

\section{Unreasonable Ubiquitousness}
\label{sec2}

We now turn to the inspiration for this paper. Three recent results exhibit quasi-polynomial behavior, in situations that seem ``unreasonable''. In particular, all three involve sets $S_t$ defined by inequalities of the form $\vec a(t)\cdot \x\le b(t)$, where $\vec a(t)$ is a polynomial in $t$; that is, the normal vectors to the facets change as $t$ changes. First we give an example showing that, unlike in Section \ref{sec1}, it is now important that we restrict to only one parameter, $t$.

\begin{example}
Define $S_{s,t}=\big\{(x,y)\in\N^2:\ sx+ty=st\big\}$. Then $S_{s,t}$ is an interval in $\Z^2$ with endpoints $(t,0)$ and $(0,s)$, and
\[\abs{S_{s,t}}=\gcd(s,t)+1.\]
\end{example}
There is no hope for simple quasi-polynomial behavior here, as the cardinality depends on the arithmetic relationship of $s$ and $t$.

\subsection{Three results}
This first result most directly generalizes Ehrhart Theory. \citet{CLS12} prove that, if $S_t$ is the set of integer points in a polytope defined by inequalities of the form $\vec a(t)\cdot \x\le b(t)$, then $\abs{S_t}$ is eventually a quasi-polynomial.

\begin{theorem}[Theorem 2.1 of \citet{CLS12}]
\label{thm:CLS}
Let $A(t)$ be an $r\times d$ matrix and $\vec b(t)$ be a column vector of length $r$, all of whose entries are in $\Z[t]$. Assume $P_t=\{\x\in\R^d:\  A(t)\x\le \vec b(t)\}$ is eventually a bounded set (a polytope). Then $\abs{P_t\cap\Z^d}$ is eventually a quasi-polynomial.
\end{theorem}

\begin{figure}
\centering
\caption{The twisting square $P_t$ from Example \ref{ex:twist}.}
\begin{subfigure}{.4\textwidth}
\centering
\includegraphics[width=.9\textwidth]{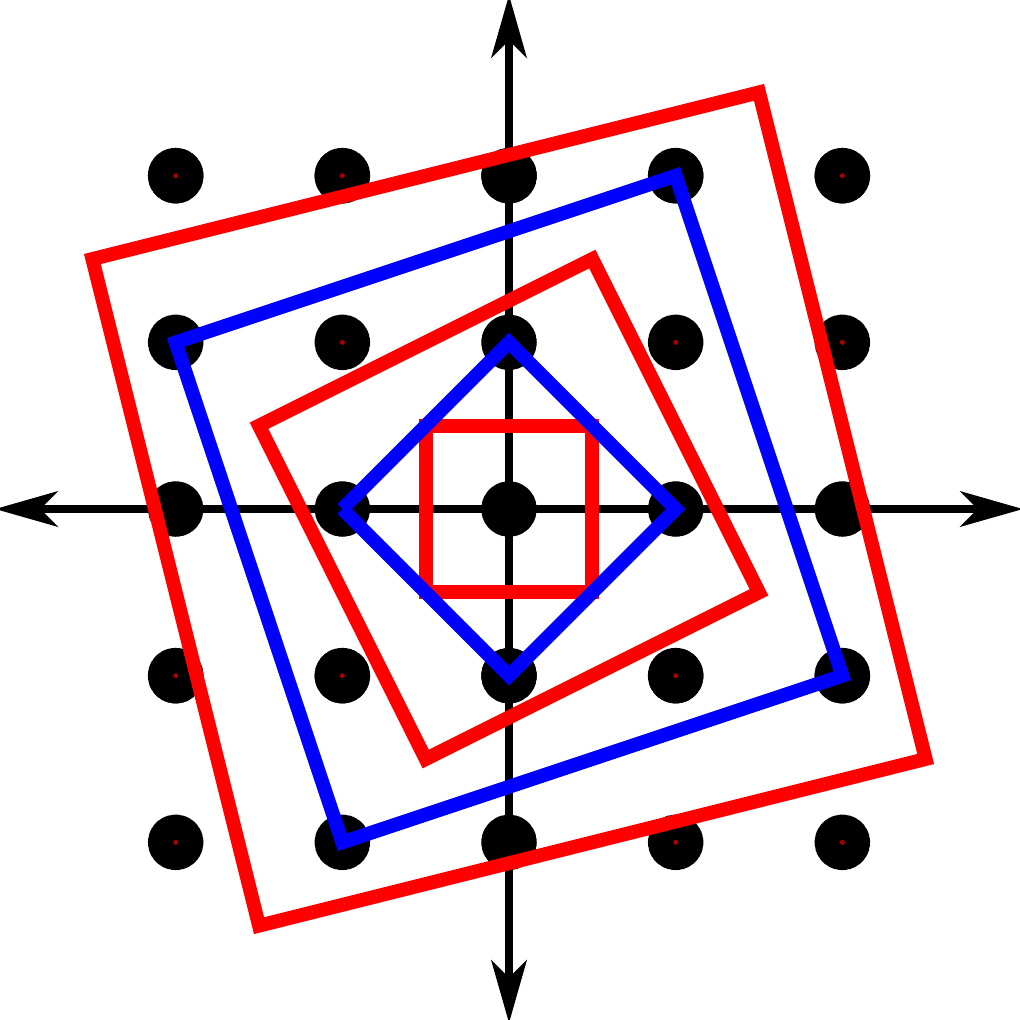}\caption{$P_t$ for $t=1,\ldots,5$.}\label{fig:twist}
\end{subfigure}
\begin{subfigure}{.4\textwidth}
\centering
\includegraphics[width=.9\textwidth]{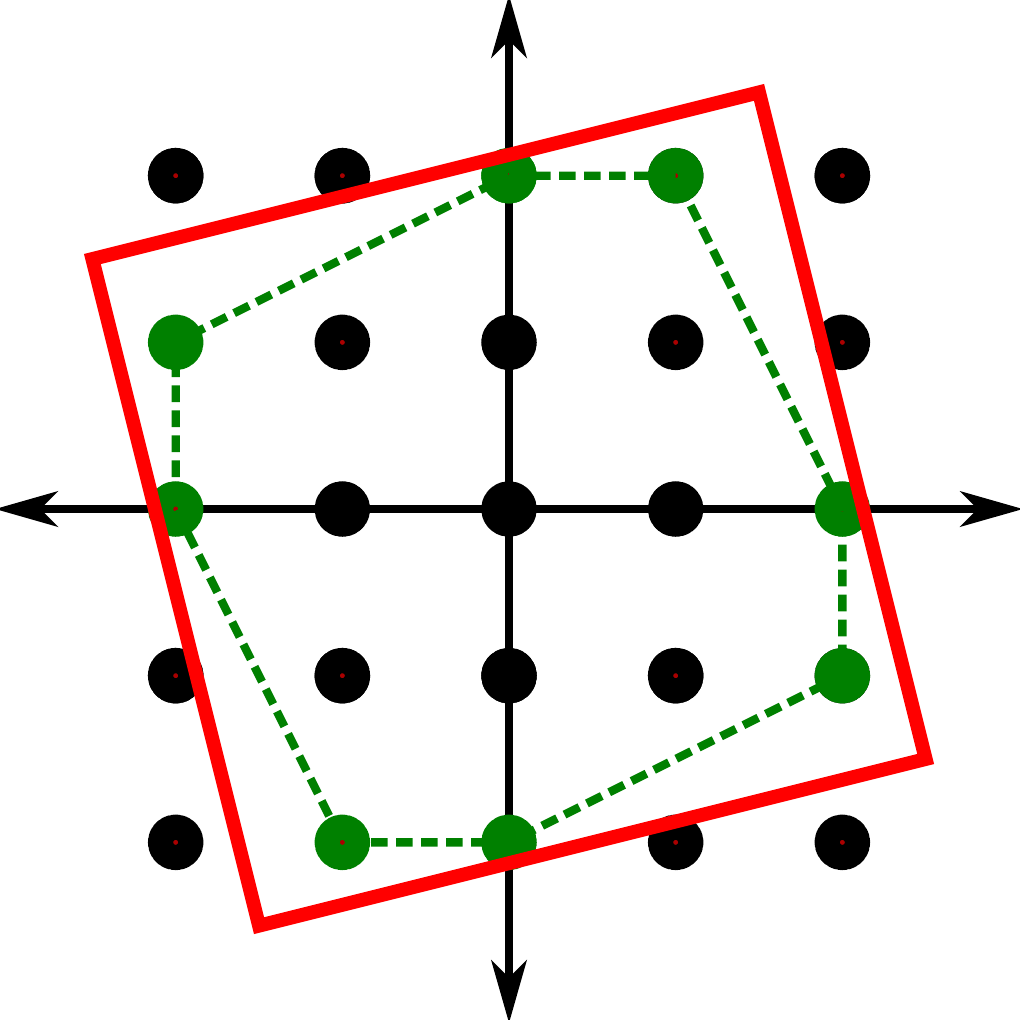}\caption{Integer hull of $P_5$.}\label{fig:twist2}
\end{subfigure}
\end{figure}

\begin{example}
\label{ex:twist}
Let $P_t$ be the ``twisting square'' in Figure \ref{fig:twist}, defined by
\[P_t=\big\{(x,y)\in\mathbb R^2:\  \abs{2x+(2t-2)y}\le t^2-2t+2,\ \abs{(2-2t)x+2y}\le t^2-2t+2\big\}.\]
Then $\abs{P_t\cap\Z^2}$ is given by the quasi-polynomial
\[\abs{P_t\cap\Z^2}=\begin{cases}t^2-2t+2 &\text{if $t$ odd,}\\t^2-2t+5 &\text{if $t$ even}. \end{cases}
\]
\end{example}

Note that Theorem \ref{thm:CLS} can be equivalently phrased  (Theorem 1.1 of \citet{CLS12}) using equalities $A(t)\x=\vec b(t)$, where $\x$ is constrained to be nonnegative, or it can be phrased (Theorem 1.4 of \citet{CLS12}) by listing the vertices of $P_t$ as rational functions of $t$.

\citet{CW11} were similarly concerned with the integer points  in polyhedra defined by $A(t)\x\le \vec b(t)$. Rather than counting $\abs{P_t\cap\Z^d}$, they wanted to find the vertices of the integer hull of $P_t$, that is, the set of vertices of the convex hull of $P_t\cap\Z^d$.

\begin{theorem}[Theorem 3.5 of \citet{CW11}]
\label{thm:CW}
Let $\vec v_i(t)$ be vectors in $\Q^d$ whose coordinates are rational functions of size O$(t)$, and let $P_t$ be the convex hull of the $\vec v_i(t)$. Then there exists a modulus $m$ and functions $\vec p_{ij}:\N \rightarrow\Z^d$ with polynomial coordinates such that, for $0\le i<m$ and for sufficiently large $t\equiv i \bmod m$, the set of vertices of the integer hull of $P_t$ is $\{\vec p_{i1}(t),\vec p_{i2}(t),\ldots,\vec p_{ik_i}(t)\}$.
\end{theorem}

\begin{example}
Consider the twisting square, $P_t$, from Example \ref{ex:twist}. When $t$ is even, the vertices of $P_t$ are integers, so the vertices of the integer hull are simply the vertices of $P_t$:
\[\left(\pm\frac{t-2}{2},\pm\frac{t}{2}\right)\quad\text{and}\quad \left(\pm\frac{t}{2},\mp\frac{t-2}{2}\right).\]
When $t$ is odd, the integer hull of $P_t$ is an octagon (pictured in Figure \ref{fig:twist2} for $t=5$) with vertices
\[\left(0,\pm\frac{t-1}{2}\right),\left(\pm\frac{t-3}{2},\pm\frac{t-1}{2}\right),\left(\pm\frac{t-1}{2},0\right),\left(\pm\frac{t-1}{2},\mp\frac{t-3}{2}\right).\]
\end{example}

Theorem \ref{thm:CW} could be similarly phrased using facet definitions of the polyhedra, rather than vertex definitions. That the vertices are O$(t)$ (grow no faster that $ct$ for some constant $c$) is important for the proof, though Calegari and Walker conjecture that the theorem still holds without this restriction.

A third recent result concerns the Frobenius number.

\begin{definition}
Given $a_1,\ldots, a_d\in\N$, let $S$ be the semigroup generated by the $a_i$, that is,
\[S=\{a\in\N:\ \exists \lambda_1,\ldots,\lambda_d\in\N,\  a=\lambda_1a_1+\cdots+\lambda_d a_d\}.\]
If the $a_i$ are relatively prime, then $S$ contains all sufficiently large integers, and the \emph{Frobenius number} is defined to be the largest integer not in $S$.
\end{definition}

Now we let $a_i=a_i(t)$ vary with $t$. \citet{RW12} prove that, if the $a_i(t)$ are linear functions of $t$, then the Frobenius number is eventually a quasi-polynomial, and they conjecture that this is true if the $a_i(t)$ are any polynomial functions of $t$:

\begin{theorem}
\label{thm:parfrob}
Let $a_i(t)\in\Z[t]$ be linear and eventually positive. Then the set of $t$ such that the $a_i(t)$ are relatively prime is eventually periodic, and, for such $t$, the Frobenius number is eventually given by a quasi-polynomial. 
\end{theorem}

\begin{example}
Consider $a_1(t)=t$, $a_2(t)=t+3$. These are relatively prime exactly when $t\equiv 1,2\bmod 3$. Since there are only two generators, a well-known formula (seemingly due to \citet{Sylvester84}) gives that the Frobenius number is
\[a_1a_2-a_1-a_2=t^2+t-3.\]
\end{example}

Note that Theorem \ref{thm:parfrob} utilizes sets defined with quantifiers; Presburger arithmetic seems a good place to look for generalizations encompassing these three results.

\subsection{Common tools}
Each of these three results has their own method for proving quasi-polynomial behavior, but there are several common tools needed. \citet{CLS12} and \citet{CW11} independently prove Theorems \ref{thm:da} through \ref{thm:dom}, \citet{CLS12} prove Theorem \ref{thm:stable}, and \citet{CW11} prove Theorem \ref{thm:round}.

\begin{theorem}[Division Algorithm] \label{thm:da}Given $f(t), g(t)$ integer-valued polynomials,
\begin{enumerate}
\item if $\deg g>0$, there exist integer-valued quasi-polynomials $q_1(t)$ and $r_1(t)$ such that $f(t)=q_1(t)g(t)+r_1(t)$, with $\deg r_1<\deg g$, and
\item if $g\ne 0$, there exist integer-valued quasi-polynomials $q_2(t)$ and $r_2(t)$ such that $f(t)=q_2(t)g(t)+r_2(t)$, with eventually $0\le r_2(t)<\abs{g(t)}$.
\end{enumerate}
\end{theorem} 

These are both useful results, and only slightly different. For example, suppose $f(t)=2t-3$ and $g(t)=t$. Then Statement 1 is a traditional polynomial division algorithm: $f=2g+ -3$. Statement 2, however, is a numerical division algorithm: $f=1g+(t-3)$, and the remainder $t-3$ is between 0 and $g$ as long as $t\ge 3$. In other words, if we have found $q_1$ and $r_1$, but we eventually have $r_1(t)<0$, then we should use quotient $q_2=q_1-\sign(g)$ and remainder $r_2=\abs{g}+r_1$ instead, as eventually $0\le \abs{g(t)}+r_1(t)<\abs{g(t)}$.

The main subtlety in proving Statement 1 of this theorem is the following: Suppose $f(t)=t^2+3t$ and $g(t)=2t+1$. Then the leading coefficient of $g$ does not divide the leading coefficient of $f$, and the traditional polynomial division algorithm would produce quotients that are not integer-valued. Instead, we look at $t$ modulo  the leading coefficient of $g$; for example, if $t$ is odd, so $t=2s+1$ for some $s\in\N$, substituting gives $f(2s+1)=4s^2+10s+3$ and $g(2s+1)=4s+3$, and now the leading term does divide evenly.

The division algorithm in hand, one can prove some stronger results:
\begin{theorem}[Euclidean Algorithm and gcds]
Let $f$ and $g$ be integer-valued quasi-polynomials. Then there exists integer-valued quasi-polynomials $p(t)$, $q(t)$, and $d(t)$ such that $\gcd\big(f(t),g(t)\big)=d(t)$  and $d(t)=p(t)f(t)+q(t)g(t)$.
\end{theorem}

This is obtained by repeated applications of the division algorithm.

\begin{example}
\[\gcd(2t+1,5t+6)=\gcd(t+4,2t+1)=\gcd(7,t+4)=\begin{cases}7&\text{ if $t\equiv 3\bmod 7$,}\\1&\text{else.}\end{cases}\]
\end{example}

Similarly, repeated application of the Euclidean algorithm can produce the Hermite or Smith normal forms of matrices. We won't define those here, but they are important, for example, in producing a basis for lower-dimensional sublattices of $\Z^d$ (see \citet{Newman72}).

\begin{theorem}[Hermite/Smith Normal Forms]
Given a matrix $A(t)$ with integer-valued quasi-polynomial entries, the Hermite and the Smith Normal forms, as well as their associated change-of-basis matrices, also have quasi-polynomial entries.
\end{theorem}

The following theorem is obvious, but is repeatedly used.
\begin{theorem}[Dominance]\label{thm:dom}
Suppose $f,g\in\Q[t]$ with $f\ne g$. Then either eventually $f(t)>g(t)$ or eventually $g(t)>f(t)$.
\end{theorem}

Repeated use of this property, for example, shows that the combinatorial structure of a polyhedron $P_t$ eventually stabilizes, when $P_t$ is defined by $A(t)\x\le \vec b(t)$:

\begin{theorem}[Stabilization]
\label{thm:stable}
Let $A(t)$ be an $r\times d$ matrix and $b(t)$ be a column vector of length $r$, all of whose entries are in $\Z[t]$, and let $P_t=\{\x\in\R^d:\  A(t)\x\le \vec b(t)\}$. For sufficiently large $t$, the set of vertices of $P_t$ are given by rational functions $\vec v_i(t)$, and the combinatorial type of $P_t$ (the subsets of vertices lying on common faces) is constant. 
\end{theorem}

Note that the vertices in this result are \emph{rational functions} of $t$: a vertex will be a point where several of the inequalities are equalities, i.e., the solution to some $A'(t)\x= \vec b'(t)$, where $A'(t)$ is a full-rank $d\times d$ matrix of polynomials in $t$. Solving for $\x$ using the adjunct matrix of $A'$ will result in $\x(t)$ given as a rational function of $t$.

For large $t$, the behavior of a rational function is predictable:

\begin{theorem}[Rounding]\label{thm:round}
Let $f(t),g(t)\in\Z[t]$. Then $f(t)/g(t)$ converges to a polynomial, and $\floor{f(t)/g(t)}$ is eventually a quasi-polynomial.
\end{theorem}

\section{Conjectures}
\label{sec3}
Let $S_t$, for $t\in\N$, be a family of subsets of $\N^d$. We now discuss some properties that it would be nice (though unreasonable!) for such sets to have, \textit{cf.} Example \ref{ex:PA}.

\begin{description}
\item[Property \hypertarget{p1}{1}:] The set of $t$ such that $S_t$ is nonempty is eventually periodic.
\end{description}
This is the weakest of the properties we will discuss, but an important one, as it is related to the decision problem -- ``Is there a solution?''

\begin{description}
\item[Property  \hypertarget{p2}{2}:] There exists a function $g:\N\rightarrow\N$ such that, if $S_t$ has finite cardinality, then $g(t)=\abs{S_t}$,  and $g(t)$ is eventually a quasi-polynomial. The set of $t$ such that $S_t$ has finite cardinality is eventually periodic.
\end{description}
This is the property found in Theorem \ref{thm:CLS}, where $S_t$ is the set of integer points in a polytope defined by inequalities of the form $\vec a(t)\cdot \x\le b(t)$. Theorems \ref{thm:CW} and \ref{thm:parfrob}, on the other hand, are not about counting points but about finding points:

\begin{description}
\item[Property  \hypertarget{p3}{3}:] There exists a function $\vec x:\N\rightarrow\N^d$ such that, if $S_t$ is non\-emp\-ty, then $\x(t)\in S_t$, and the coordinate functions of $\x$ are eventually quasi-polynomials. The set of $t$ such that $S_t$ is nonempty is eventually periodic.
\end{description}
This function $\x(t)$ acts as a certificate that the set is nonempty. But we may want to go further and pick out particular elements of $S_t$:

\begin{description}
\item[Property  \hypertarget{p3p}{3a}:] Given $\vec c\in \Z^d\setminus\{0\}$, there exists a function $\vec x:\N\rightarrow\N^d$ such that, if $\max_{\y\in S_t} \vec c\cdot \y$ exists, then it is attained at $\x(t)\in S_t$, and the coordinate functions of $\x$ are eventually quasi-polynomials. The set of  $t$ such that the maximum exists is eventually periodic.
\end{description}
This corresponds to Theorem \ref{thm:parfrob}, where we want to find the Frobenius number, the maximum element of the complement of the semigroup. On the other hand, we may want to list \emph{multiple} elements of the set:

\begin{description}
\item[Property  \hypertarget{p3pp}{3b}:] Fix $k\in\N$. There exist functions $\vec x_1,\ldots,\vec x_k:\N\rightarrow\N^d$ such that, if $\abs{S_t}\ge k$, then $\x_1(t),\ldots,\x_k(t)$ are distinct elements of $S_t$, and the coordinate functions of $\x_i$ are eventually quasi-polynomials. The set of $t$ such that $\abs{S_t}\ge k$ is eventually periodic.
\end{description}
If there is a uniform bound on $\abs{S_t}$, then this property can be used to enumerate all elements of $S_t$, for all $t$. This is the content of Theorem \ref{thm:CW}. Property 2 is about counting all solutions and Properties 3/3a/3b are about obtaining specific solutions, and so they seem somewhat orthogonal to each other. The following property, we shall see, unifies them:

\begin{description}
\item[Property  \hypertarget{p4}{4}:] There exists a period $m$ such that, for $t\equiv i\bmod m$, \[\sum_{\vec s\in S_t}\x^{\vec s} = \frac{\sum_{j=1}^{n_i}\alpha_{ij}\x^{\vec q_{ij}(t)}}{(1-\vec x^{\vec b_{i1}(t)})\cdots(1-\vec x^{\vec b_{ik_i}(t)})},\]
where $\alpha_{ij}\in\Q$, and the coordinate functions of $\vec q_{ij},\vec b_{ij}:\N\rightarrow\Z^d$ are polynomials with the $\vec b_{ij}(t)$ lexicographically positive.
\end{description}

For what sets $S_t$ can we hope for these properties to hold? Here is a candidate:

\begin{definition}
A family of sets $S_t$ is a \emph{parametric Presburger family} if they can be defined over the natural numbers using quantifiers, boolean operations, and  inequalities of the form $\vec a(t)\cdot \x\le b(t)$, where $b\in\Z[t]$ and $\vec a\in\Z[t]^d$.
\end{definition}

We conjecture that these properties do, in fact, hold for any parametric Presburger family:

\begin{conjecture}
Let $S_t$ be a parametric Presburger family. Then Properties 1, 2, 3, 3a, 3b, and 4 all hold.
\end{conjecture}

Note that one can define a family $S_t$ of subsets of $\Z^d$ rather than of $\N^d$, but care must be taken with generating functions; see Remark \ref{rmk:lines}.

As evidence that Property 4 is interesting, we will show that it generalizes both 2 and 3/3a/3b:

\begin{theorem}
\label{thm:implications}
Let $S_t$ be any family of subsets of $\N^d$. We have the following implications among possible properties of $S_t$.
\[
\xymatrix@R.1in@C.1in{ & & & \text{2} \ar@{=>}@/^/[drrr]  & & &\\
\text{4} \ar@{=>}@/_/[drr] \ar@{=>}[rr] \ar@{=>}@/^/[urrr] & & \text{3a} \ar@{=>}[rr] &  & \text{3} \ar@{=>}[rr]  & & \text{1}\\
& & \text{3b} \ar@{=>}@/_/[urr] &  & & }
\]
\end{theorem}

As a final relationship between these properties, we note that, for the class of parametric Presburger families, 3, 3a, and 3b are equivalent:

\begin{theorem}
\label{thm:equiv}
Suppose all parametric Presburger families have Property 3. Then all parametric Presburger families have Properties 3a and 3b.
\end{theorem}

Theorem \ref{thm:equiv} is a weaker implication than Theorem \ref{thm:implications}, which holds for a single family $S_t$ in isolation. To prove that 3 ``implies'' 3a and 3b, on the other hand, we will need to create new families ${S}'_t$ using additional quantifiers or boolean operators, and we need to know that these new families still have Property 3.

Finally, we give evidence that  these properties might actually hold. We can show that they all hold for two broad classes of parametric Presburger families:

\begin{theorem}
\label{thm:cases}
Suppose $S_t$ is a parametric Presburger family such that either
\begin{enumerate}
\item[(a)] $S_t$ is defined without using any quantifiers, \emph{or}
\item[(b)] the only inequalities used to define $S_t$ are of the form $\vec a\cdot \x\le b(t)$, where $b(t)$ is a polynomial (that is, the normal vectors to the hyperplanes must be fixed).
\end{enumerate}
Then Properties 1, 2, 3, 3a, 3b, and 4 all hold.
\end{theorem}

We isolate a piece of the proof of Part (b), in order to point out that Property 4 is a weaker property than we might hope for, but seems to be as strong a property as we can get. Indeed, we might hope that $\sum_{t\in\N, \vec s\in S_t}\x^{\vec s}y^t$ is a rational generating function. Theorem \ref{thm:PA} shows that this is true for sets defined in the normal Presburger arithmetic, and the following theorem shows that this implies Property 4.

\begin{theorem}
\label{thm:diff_gfs}
Suppose $S_{\vec p}$, for $\vec p\in\N^n$, is a family of subsets of $\N^d$. If $\sum_{\vec p\in\N^n, \vec s\in S_{\vec p}}\x^{\vec s}\y^{\vec p}$ is a rational generating function, then there is a finite decomposition of $\N^n$ into pieces of the form $P\cap\Z^n$ (with $P$ a polyhedron) such that, considering the $\vec p$ in each piece separately,
\[\sum_{\vec s\in S_{\vec p}}\x^{\vec s} = \sum_{i=1}^m\epsilon_i\frac{\x^{\vec q_i(\vec p)}}{(1-\vec x^{\vec b_{i1}})\cdots(1-\vec x^{\vec b_{ik_i}})},\]
where $\epsilon_i=\pm1$, $\vec b_{ij}\in\Z^d$ are lexicographically positive, and the coordinate functions of $\vec q_i:\N^n\rightarrow\Z^d$ are linear quasi-polynomials in $\vec p$.
\end{theorem}

In general, however, $\sum_{t\in\N, \vec s\in S_t}\x^{\vec s}y^t$ will not be a rational generating function:

\begin{example}
Let $S_t$ be the set $\{(s_1,s_2)\in\N^2:\ ts_1=s_2\}$. Then
\[\sum_{\vec s\in S_t}\x^{\vec s}=1+x_1x_2^t+x_1^2x_2^{2t}+\cdots=\frac{1}{1-x_1x_2^t}\]
is a rational generating function with exponents depending on $t$, so Property 4 is satisfied. Nevertheless,
\[\sum_{t\in\N,\vec s\in S_t} \x^{\vec s}y^t=\frac{1}{1-x_1}+\frac{y}{1-x_1x_2}+\frac{y^2}{1-x_1x_2^2}+\cdots\]
cannot be written as a rational function. To prove that it cannot be so written, note that the set $\big\{(s_1,s_2,t):\ \vec s\in S_t\big\}$ cannot be written as a finite union of sets of the form $P\cap(\vec\lambda+\Lambda)$, where $P$ is a polyhedron, $\vec\lambda\in\Z^3$ and $\Lambda\subseteq\Z^3$ is a lattice; Theorem 1 of \citet{Woods13} then implies that $\sum_{t\in\N,\vec s\in S_t}\x^{\vec s}y^t$ is not a rational generating function.
\end{example}

\section{Proofs}
\label{sec4}
Because of dependencies among these proofs, we present them out of numerical order. We first prove Theorem \ref{thm:diff_gfs}. Then we prove that Property 4 holds for the sets given in Theorem \ref{thm:cases}. Next we prove Theorem \ref{thm:implications}. The rest of the properties in Theorem \ref{thm:cases} are an immediate corollary. Finally, we prove Theorem \ref{thm:equiv}.

\subsection{Proof of Theorem \ref{thm:diff_gfs}}\label{pf3.6}
We want to compute the coefficient of $\y^{\vec p}$ as a rational generating function in $\x$. We can reduce to computing the coefficient of $\y^{\vec p}$ for a single term of the form
\begin{equation}
\label{eqn:gf_allp}
\frac{\x^{\vec q}\y^{\vec r}}{(1-\x^{\vec c_1}\y^{\vec d_1})\cdots (1-\x^{\vec c_k}\y^{\vec d_k})},
\end{equation}
where $\vec q,\vec c_i\in\N^d$, $\vec r,\vec d_i\in\N^n$, and $\vec c_i,\vec d_i$ are nonzero: indeed, by the assumption that $\sum_{\vec p\in\N^n, \vec s\in S_{\vec p}}\x^{\vec s}\y^{\vec p}$ is a rational generating function, it can be written as a linear combination of such terms as (\ref{eqn:gf_allp}), so we could examine each such term separately and then apply the linear combination.

Let $\vec p\in\N^n$ be fixed. Expanding (\ref{eqn:gf_allp}) as a product of geometric series, we see that we get a term $\x^{\vec s}\vec y^{\vec p}$ exactly when there exist $\lambda_i\in \N$ such that
\begin{equation}\label{eqn:1term}
\vec p = \vec r+ \sum_{i}\lambda_i\vec d_i\quad\text{and}\quad \vec s=\vec q + \sum_i\lambda_i\vec c_i.
\end{equation}
Define the parametric polyhedron
\[Q_{\vec p}=\Big\{\vec \lambda\in\R^k:\ \lambda_i\ge 0\text{ and }\vec p = \vec r+ \sum_{i}\lambda_i\vec d_i\Big\}.\]
Then (\ref{eqn:1term}) tells us that, for every $\vec\lambda\in Q_{\vec p}\cap\Z^k$, we will get a monomial $\x^{\vec q+\sum_i\lambda_i\vec c_i}\y^{\vec p}$ in the expansion of (\ref{eqn:gf_allp}), and therefore
\[\sum_{\vec s\in S_{\vec p}}\x^{\vec s}=\sum_{\vec \lambda\in Q_{\vec p}\cap\Z^k}\x^{\vec q+\sum_i\lambda_i\vec c_i}.\] Theorem \ref{thm:VW} implies that there is a finite decomposition of $\Z^n$ into pieces of the form $P\cap\Z^n$ (with $P$ a polyhedron) such that, considering the $\vec p$ in each piece separately,
\[\sum_{\vec{\lambda}\in Q_{\vec p}\cap\Z^k}\z^{\vec \lambda}=\sum_{i=1}^m\epsilon_i\frac{\vec{ z}^{\vec u_i(\vec p)}}{(1-\vec{z}^{\vec b_{i1}})\cdots(1-\vec{z}^{\vec b_{ik_i}})},\]
where $\epsilon_i=\pm1$, the coordinates of $\vec u_i$ are linear quasi-polynomials in $\vec p$, and $\vec b_{ij}\in\Z^d$ are lexicographically positive. Substituting $z_i=\vec x^{\vec c_i}$ and multiplying by $\vec x^{\vec q}$ yields $\sum_{\lambda\in Q_{\vec p}\cap\Z^k}\x^{\vec q+\sum_i\lambda_i\vec c_i}.$ This is the desired generating function, $\sum_{\vec s\in S_{\vec p}}\x^{\vec s}$, in the form required by the theorem.

\subsection{Proof of Property 4 in Theorem \ref{thm:cases}}\label{pf3.54}

\noindent {\bf Part (a):}
We want to apply a key idea of \citet[Lemma 3.2]{CLS12}, but first we need to simplify by reducing our problem (that of showing Property 4 holds for any quantifier-free Presburger formula) to simpler ones.

Let the set of inequalities used in the quantifier-free Presburger formula be given by $\{\vec a_i(t)\cdot \x\le b_i(t)\}$, $1\le i\le n$. Our domain $\N^d$ is partitioned into (at most) $3^n$ pieces, as follows: to create a piece, for each $i$, choose either $\vec a_i(t)\cdot \x< b_i(t)$, $\vec a_i(t)\cdot \x> b_i(t)$, or $\vec a_i(t)\cdot \x= b_i(t)$, and let a piece of the  partition be the set of $\x$ satisfying the conjunction of these (in)equalities. Within each piece, either all points will satisfy the Presburger formula or all will fail to satisfy the Presburger formula. It suffices to prove Property 4 for one of these pieces that satisfy the Presburger formula, because the final generating function is simply a sum of the generating functions of these pieces. These pieces are polyhedra (since we are looking at integer points, we may replace the open 
$\vec a_i(t)\cdot \x< b_i(t)$ with the closed $\vec a_i(t)\cdot \x\le b_i(t)-1$), so it suffices to prove Property 4 for polyhedra, $P_t$, defined with equations of the form $\vec a(t)\cdot \x\le b(t)$.

We still can't apply \citet{CLS12}, because that will only apply to bounded polyhedra. For sufficiently large $t$, Theorem \ref{thm:stable} shows that the combinatorics of the vertex cones of $P_t$ stabilize.  Since Theorem \ref{thm:Brion} (Brion's Theorem) gives that the generating function for $P_t$ is the sum of the generating functions for these vertex cones, it suffices to prove Property 4 for such cones
\[K_t=\vec v(t)+\cone\big(\vec u_1(t),\ldots,\vec u_n(t)\big)=\big\{\vec v(t)+\sum_i\lambda_i\vec u_i(t):\ \lambda_i\ge 0\big\}\]
where the coordinates of $\vec v$ are rational functions and the coordinates of $\vec u_i$ are polynomials.

Now take $\Pi_t$ to be the fundamental parallelepiped of $K_t$, that is,
\[\Pi_t=\big\{\vec v(t)+\sum_i\lambda_i\vec u_i(t):\  0\le\lambda_i<1\big\},\]
and note that
\[\sum_{\vec s\in K_t\cap\Z^d}\vec x^{\vec s}=\frac{\sum_{\vec s\in \Pi_t\cap\Z^d}\vec x^{\vec s}}{\left(1-\vec x^{\vec u_1(t)}\right)\cdots \left(1-\vec x^{\vec u_n(t)}\right)}.\]
It therefore suffices to prove Property 4 for $\Pi_t$, which is a \emph{bounded} polyhedron, so \citet{CLS12} finally applies, as follows. Given $k$, define the map $\vec\psi:\Z^{d(k+1)}\rightarrow\Z^d$, taking $\bar{\vec s}=(\bar{s}_{ij})_{1\le i\le d,0\le j\le k}$ to $\vec s=(s_i)_{1\le i\le d}$, given by
\[\big(\vec\psi(\bar{\vec s})\big)_i=\bar{s}_{ik}t^k+\cdots+\bar{s}_{i1}t+\bar{s}_{i0}.\]
Then Lemma 3.2 of \citet{CLS12}  states that, for some $k$ and some $R_t\subseteq\Z^{d(k+1)}$, $\vec\psi:R_t\rightarrow\Pi_t\cap\Z^d$ is a bijection, and furthermore that $R_t$ can be defined as the disjoint union of the integer points in polyhedra which satisfy particular linear inequalities of the form
\[\bar{\vec a}\cdot \bar{\vec x}\le b(t),\]
where $\bar{\vec x}=(\bar{x}_{ij})_{ij}$, $b(t)$ is \emph{linear}, and $\bar{\vec a}$ no longer depends on $t$. Then Theorem \ref{thm:VW} shows that (for sufficiently large $t$ so that the combinatorial structure of this set stabilizes)
\begin{equation}\label{eqn:CLS}
\sum_{\bar{\vec{s}}\in R_t}\vec{\bar{x}}^{\bar{\vec{s}}}=\sum_{i=1}^m\epsilon_i\frac{\vec{\bar x}^{\vec u_i(t)}}{(1-\vec{\bar{x}}^{\vec d_{i1}})\cdots(1-\vec{\bar{x}}^{\vec d_{ik_i}})},\end{equation}
where $\epsilon=\pm1$, the coordinates of $\vec u_i$ are linear quasi-polynomials in $t$, and $\vec d_{ij}\in\Z^{d(k+1)}$ are nonzero.
Substituting $\bar{x}_{ij}=x_i^{t^j}$ into a monomial $\bar{\vec x}^{\bar{\vec s}}$ gives
\[\prod_{i,j}x_i^{\bar s_{ij}t^j}\prod_ix_i^{\bar{s}_{ik}t^k+\cdots+\bar{s}_{i1}t+\bar{s}_{i0}}=\vec x^{\vec \psi(\bar{\vec s})}.\]
Therefore substituting $\bar{x}_{ij}=x_i^{t^j}$ into (\ref{eqn:CLS}) shows that $\sum_{\vec s\in \Pi_t\cap\Z^d}\vec x^{\vec s}$ has Property 4, as desired.
\vspace{.1in}

\noindent {\bf Part (b):}  Let the set of inequalities used in the Presburger formula be given by $\{\vec a_i\cdot \x\le b_i(t)\}$, $1\le i\le n$. For the moment, replace each $b_i(t)$ with a single variable $p_i$, so that we have a Presburger formula $F(\vec x, \vec p)$, defined using inequalities $\vec a_i\cdot x\le p_i$. Let
\[S_{\vec p}=\big\{\x\in\N^d:\ F(\vec x, \vec p)\big\},\]
so that $S_t=S_{(b_1(t),\ldots,b_n(t))}$.

Theorem \ref{thm:PA} implies that $\sum_{\vec p\in\N^n,\vec s\in S_{\vec p}}\x^{\vec s}\y^{\vec p}$ is a rational generating function. Theorem \ref{thm:diff_gfs} then implies that there is a finite decomposition of $\N^n$ into pieces of the form $P\cap\Z^n$ (with $P$ a polyhedron), such that, considering the $\vec p$ in each piece separately,
\[\sum_{\vec s\in S_{\vec p}}\x^{\vec s} = \sum_{i=1}^m\epsilon_i\frac{\x^{\vec q_i(\vec p)}}{(1-\vec x^{\vec d_{i1}})\cdots(1-\vec x^{\vec d_{ik_i}})},\]
where $\epsilon_i=\pm1$, $\vec d_{ij}\in\Z^d$ are lexicographically positive, and the coordinate functions of $\vec q_i$ are linear quasi-polynomials in $\vec p$. For sufficiently large $t$, the vector $\big(b_1(t),\ldots,b_n(t)\big)$ eventually stays in one of these pieces of the decomposition. Therefore, we may substitute $\vec p=\big(b_1(t),\ldots,b_n(t)\big)$ into the appropriate generating function. This proves Property \hyperlink{4}{4} for $S_t$.

\subsection{Proof of Theorem \ref{thm:implications}}\label{pf3.3}
\noindent \textbf{\hyperlink{p2}{2} $\mathbf{\Rightarrow}$ \hyperlink{p1}{1}} : 
Suppose $g(t)=\abs{S_t}$ is eventually a quasi-polynomial. Let $m$ be a period and $p_i(t)$ be polynomials such that, eventually, $g(t)=p_i(t)$ for $t\equiv i \bmod m$. Let $I$ be the indices $i$ such that $p_i$ is not identically 0. For $i\in I $, we have $\abs{S_t}=p_i(t)>0$ for sufficiently large $t$, so, eventually, $S_t$ is nonempty if and only if $t \bmod m$ is in $I$. 
\vspace{.1in}

\noindent \textbf{\hyperlink{p3}{3} $\mathbf{\Rightarrow}$ \hyperlink{p1}{1}} : 
 Property 1 is in the definition of Property 3.

\vspace{.1in}

\noindent \textbf{\hyperlink{p3p}{3a} $\mathbf{\Rightarrow}$ \hyperlink{p3}{3}} : 
 Take $\vec c=(-1,\ldots,-1)$. Then $S_t\subseteq\N^d$ is nonempty if and only if $\vec c\cdot \y$ has a maximum on $S_t$. By property 3a, the set of such $t$ is eventually periodic, and when nonempty, 3a gives an element $\x(t)\in S_t$.

\vspace{.1in}

\noindent \textbf{\hyperlink{p3pp}{3b} $\mathbf{\Rightarrow}$ \hyperlink{p3}{3}} : 
 Property 3 is Property 3b with $k=1$.

\vspace{.1in}

\noindent \textbf{\hyperlink{p4}{4} $\mathbf{\Rightarrow}$ \hyperlink{p2}{2}} : 
 Let $f(\x)=\sum_{\vec s\in S_t}\x^{\vec s}$. Let us examine a particular residue class, so that, for sufficiently large $t\equiv i\bmod m$, we have
\[f(\x) = \frac{\sum_j\alpha_j\x^{\vec q_j(t)}}{(1-\vec x^{\vec b_1(t)})\cdots(1-\vec x^{\vec b_k(t)})},\]
where $\alpha_j\in\Q$, and the coordinate functions of $\vec q_j,\vec b_j:\N\rightarrow\N^d$ are polynomials. When $f$ is defined at $(1,\ldots,1)$,  $f(1,\ldots,1)=\abs{S_t}<\infty$, and when $f$ has a pole at $(1,\ldots,1)$, $S_t$ has infinite cardinality. So we attempt to compute $\lim_{\x\rightarrow \vec 1}f(\x)$, letting $x_j\rightarrow 1$ one variable at a time, repeatedly using L'H\^opital's rule.

We will repeatedly take partial derivatives of the numerator (and denominator) and set some of the $x_j$ to 1. At each point, we will see that the numerator (and denominator) will look like
\begin{equation}\label{eqn:plugging_in}g(\bar{\x})=\sum_j \beta_j(t)\bar{\x}^{\vec r_j(t)},\end{equation}
where $\beta_j$ and the coordinate functions of $\vec r_j$ are polynomials, and $\bar\x$ is a subset of the variables of $\x$. Our base case, $f$, certainly has numerator and denominator of this form. 

Note that, at each step, we may simplify (\ref{eqn:plugging_in}) so that the $\vec r_j$ are distinct polynomials and the $\beta_j$ are not identically zero. Then for sufficiently large $t$, $\bar\x^{\vec r_j(t)}$ are distinct monomials and $\beta_j(t)$ are nonzero. Therefore $g(\bar\x)$ is either a nonzero polynomial in $\bar\x$ for all sufficiently large $t$ or identically zero for all sufficiently large $t$.

Starting with $x_1$, we substitute $x_1=1$ into the numerator and denominator of $f$. If both are eventually nonzero, we continue to $x_2$. If the numerator is eventually nonzero but the denominator is eventually zero, then $(1,\ldots,1)$ is eventually a pole of $f(\x)$, meaning that for sufficiently large $t\equiv i\bmod m$, $\abs{S_t}$ is infinite. If the numerator is eventually zero but the denominator is eventually nonzero, then $\abs{S_t}=f(1,\ldots,1)=0$. If both the numerator and denominator are eventually 0, then we apply L'H\^opital's rule, taking the derivatives of  the numerator and denominator with respect to $x_1$. Our new numerator and denominator will continue having the form (\ref{eqn:plugging_in}).

We continue with $x_1$ in this way, and then continue with $x_2$, etc., until we have either established that $\abs{S_t}$ is infinite for all $t\equiv i\bmod m$, zero for all $t\equiv i\bmod m$, or we have eliminated all of the variables among $\x$. In this last case, the numerator and the denominator of the result are still in the form  (\ref{eqn:plugging_in}), so we must have $\abs{S_t}=g(t)/h(t)$, where $g,h\in\Q[t]$. To conclude, we need to show that $g(t)/h(t)$ must actually be a polynomial. Indeed, this follows from the fact that $g(t)/h(t)$ must always be an integer, $\abs{S_t}$: Apply the standard polynomial division algorithm in $\Q[t]$ so that
\[\frac{g(t)}{h(t)}=q(t)+\frac{r(t)}{h(t)},\]
with $\deg r<\deg h$. Let $n$ be a common multiple of the denominators of the coefficients of $q(t)$, so that $nq(t)$ is always an integer. Since $ng(t)/h(t)$ is always an integer, $nr(t)/h(t)$ must also always be an integer. Since $\deg r<\deg h$, eventually $\abs{nr(t)/h(t)}<1$, so the only possibility is that $r=0$. Then $\abs{S_t}=g(t)/h(t)=q(t)$ is a polynomial (for $t\equiv i\bmod m$), and we have proven Property \hyperlink{p2}{2}.

\vspace{.1in}

\noindent \textbf{\hyperlink{p4}{4} $\mathbf{\Rightarrow}$ \hyperlink{p3p}{3a}} :

Without loss of generality, assume $c_1\ne 0$ (otherwise, reorder the variables). Define the map $\vec\phi:\N^d\rightarrow\Z\times\N^{d-1}$ by
\[\vec\phi(\vec x)=(-\vec c\cdot\vec x,x_2,x_3,\ldots,x_d).\]
Since $c_1\ne 0$, this map is injective, and therefore the map is bijective from $S_t$ onto its image $\vec\phi(S_t)$.

Assume that $\vec c\cdot \vec x$ has a maximum on $S_t$. Then $\vec\phi(S_t)$ has a unique lexicographically minimal element; this is the element we will give a function for. Given $\sum_{\vec s\in S_t}\vec x^{\vec s}$ as in Property 4, substitute in $x_1=y_1^{-c_1}$ and $x_i=y_1^{-c_i}y_i$ for $2\le i\le d$; note that
\[x_1^{s_1}\cdots x_d^{s_d}\quad\text{becomes}\quad y_1^{-c_1s_1}\left(y_1^{-c_2s_2}y_2^{s_2}\right)\cdots\left(y_1^{-c_ds_d}y_d^{s_d}\right)=y_1^{-\vec c\cdot \vec s}y_2^{s_2}\cdots y_d^{s_d}=\y^{\vec\phi(\vec s)},\]
and so the generating function becomes
\begin{equation}\label{eqn:lexmin}\sum_{\vec p\in\vec\phi(S_t)} \vec y^{\vec p} = \frac{\sum_{j=1}^{n}\alpha_{j}\y^{\vec r_{j}(t)}}{(1-\vec y^{\vec d_{1}(t)})\cdots(1-\vec y^{\vec d_{k}(t)})}.\end{equation}
Restricting our attention to sufficiently large $t$ and a particular residue class of $t\bmod m$, we have $\alpha_{j}\in\Q\setminus\{0\}$ and the coordinate functions of $\vec r_{j},\ \vec d_{j}$, are polynomials. We may assume that $\vec d_{j}(t)$ are lexicographically positive, by taking sufficiently large $t$ so that the signs of the coordinate functions do not change, and then using Remark \ref{rmk:lexpos}.

Since we are assuming that $\vec c\cdot\vec x$ has a maximum, $M$, on $S_t$, we see that all $\vec p\in\vec\phi(S_t)$ have $(-M,0,0,\ldots,0)\le \vec p$ coordinatewise. Therefore $\sum_{\vec p\in\vec\phi(S_t)} \vec y^{\vec p}$ is a Laurent power series convergent on a neighborhood of  $\vec y = (e^{-\varepsilon},e^{-\varepsilon^2},\ldots,e^{-\varepsilon^d})$, for any $\varepsilon>0$. Expanding the right-hand side of (\ref{eqn:lexmin}) as a product of geometric series also gives a Laurent power series convergent on a neighborhood of  $\vec y = (e^{-\varepsilon},e^{-\varepsilon^2},\ldots,e^{-\varepsilon^d})$, for sufficiently small $\varepsilon$ (since the $\vec d_j$ are lexicographically positive; see Remark \ref{rmk:lexpos2}). Therefore they must be identical as Laurent power series expansions. But the lexicographically minimal term of the expanded right-hand side is clear: for sufficiently large $t$, it is the unique lexicographical minimum among the $\vec r_j(t)$, $1\le j\le n$. Inverting $\vec\phi$ yields an element of $S_t$ maximizing $\vec c \cdot \vec x$.

Finally, we will show that, $\vec c \cdot \vec x$ either eventually always achieves a maximum or eventually always doesn't (at least for this residue class of $t\bmod m$; this implies that the set of \emph{all} $t$ for which it achieves a maximum is eventually periodic). If the maximum exists, we have found it above; call it $M(t)$, which is a polynomial. Define the set $R_t=\{\vec x\in\N^d:\ \vec c\cdot \vec x \ge M(t)+1\}$. Then $\vec c\cdot \vec x$ attains a maximum on $S_t$ if and only if $S_t\cap R_t$ is empty. 

Given two Laurent power series $f(\vec s)=\sum_{\vec s\in\Z^d} p_{\vec s}\x^{\vec s}$ and $g(\vec x)=\sum_{\vec s\in\Z^d} q_{\vec s}\x^{\vec s}$ ($p_{\vec s}, q_{\vec s}\in \Q$) convergent on a neighborhood of some $(e^{-\varepsilon},e^{-\varepsilon^2},\ldots,e^{-\varepsilon^d})$, define their Hadamard product, $f\star g$, to be $\sum_{\vec s\in\Z^d} p_{\vec s}q_{\vec s}\x^{\vec s}$. We are given that $f(\x)\doteq\sum_{\vec s\in S_t}\vec x^{\vec s}$ has Property 4. By Theorem \ref{thm:VW}, $g(\x)\doteq\sum_{\vec s\in R_t}\vec x^{\vec s}$ also has Property 4. In particular, their numerator and denominator have the form
\begin{equation}\label{eqn:plugging_in2}\sum_j \beta_j(t)\x^{\vec r_j(t)},\end{equation}
where $\beta_j$ and the coordinate functions of $\vec r_j$ are polynomials.  Note that $f\star g=\sum_{\vec s\in S_t\cap R_t}\vec x^{\vec s}$; we want to decide whether this is identically zero. We will show that $f\star g$ also has numerator and denominator in the form (\ref{eqn:plugging_in2}), in which case (for sufficiently large $t$) it is identically zero if and only if its numerator is identically zero.

Since the Hadamard product is a bilinear operator, it suffices to show that
\[\frac{\x^{\vec q(t)}}{(1-\x^{\vec b_1(t)})\cdots (1-\x^{\vec b_{k}(t)})}\star \frac{\x^{\vec r(t)}}{(1-\x^{\vec d_1(t)})\cdots (1-\x^{\vec d_\ell(t)})}\]
has numerator and denominator of the form (\ref{eqn:plugging_in2}). What is the coefficient of $\x^{\vec s}$ in this Hadamard product? Expanding the first rational function as a product of geometric series, the coefficient of $\vec x^{\vec s}$ in it is
\[\#\vec \lambda\in\N^k:\ \vec q(t)+\lambda_1\vec b_1(t)+\cdots+\lambda_k\vec b_k(t)=\vec s,\]
and similarly the coefficient of  $\vec x^{\vec s}$ in the second is $\#\vec\mu\in\N^\ell:\ \vec r+\vec\mu\cdot\vec d(t)=\vec s$. Therefore the coefficient of $\vec x^{\vec s}$ in the Hadamard product is the product of these, that is,
\begin{equation}\label{eqn:had}\#(\vec\lambda,\vec\mu)\in\N^k\times\N^\ell:\ \vec q+\vec\lambda\cdot\vec b(t)=\vec r+\vec\mu\cdot\vec d(t)=\vec s.\end{equation}

There are no quantifiers in this expression (\ref{eqn:had}), so it is of the form in Part (a) of Theorem \ref{thm:cases}, meaning we proved in Section \ref{pf3.54} that the generating function
\[\sum_{\vec s\in\N^d,\vec\lambda\in\N^k,\vec\mu\in\N^\ell\text{ satisfying (\ref{eqn:had})}}\x^{\vec s}\y^{\vec \lambda}\z^{\vec\mu}\]
satisfies Property 4. Substituting $\y=\vec 1,\z=\vec 1$ gives the Hadamard product; following the proof of $4\Rightarrow 2$ (which may involve L'H\^opital's rule and therefore differentiation of the numerator and denominator), we see that the Hadamard product indeed has numerator and denominator in the form (\ref{eqn:plugging_in2}). Therefore $f\star g$ is either identically zero (and hence $\vec c\cdot\x$ is always maximized on $S_t$) or eventually never zero (and hence $\vec c\cdot\x$ does not achieve a maximum). This is true for every residue class of $t$, so, all together, the set of $t$ for which a maximum is achieved will eventually be periodic.

\vspace{.1in}

\noindent \textbf{\hyperlink{p4}{4} $\mathbf{\Rightarrow}$ \hyperlink{p3pp}{3b}} : 
Since we have proven that \hyperlink{4}{4} $\Rightarrow$ \hyperlink{p3p}{3a} $\Rightarrow$ \hyperlink{3}{3}, we can get a single element $\y_1(t)\in S_t$, when $S_t$ is nonempty. Now we can continue inductively: $\sum_{\vec s\in S_t\setminus\{\y_1(t)\}}\x^{\vec s}=\left(\sum_{\vec s\in S_t}\x^{\vec s}\right) - \x^{\y_1(t)}$ can be written as a rational function of the appropriate form, so $S_t\setminus\{\y_1(t)\}$ has Property 4. Thus we can use Property 3 to get another element $\y_2(t)$ and also conclude that the set of $t$ for which $\abs{S_t}\ge 2$ (which is the set of $t$ for which $\abs{S_t\setminus\{\y_1(t)\}}$ is nonempty) is eventually periodic. Continue inductively.

\subsection{Proof of Theorem \ref{thm:cases}} \label{pf3.5}
Section \ref{pf3.54} proved Property 4 for the sets given in Theorem \ref{thm:cases}. The other properties follow immediately using Theorem \ref{thm:implications}.

\subsection{Proof of Theorem \ref{thm:equiv}} \label{pf3.4}
Assume that Property \hyperlink{p3}{3} holds for all parametric Presburger families. For $t\in\N$, let $S_t=\{\x\in\N^d:\ F(\x,t)\}$ be a parametric Presburger family defined by formula $F$. We must show that \hyperlink{p3p}{3a} and \hyperlink{p3pp}{3b} hold for $S_t$.

To prove that 3a holds, let $\vec c\in\Z^d\setminus\{0\}$ be given. Define a new set ${S}'_t$ to be the solutions $\vec x$ to
\[F(\x,t)\wedge\forall \vec y\big(F(\vec y,t) \Rightarrow \vec c\cdot \vec y \le \vec c \cdot \x\big).\]
These are exactly the set of elements $\x$ maximizing $\vec c\cdot \x$ in $S_t$ (and it is empty if there is no such $\x$). By Property {3}, the set of $t$ such that ${S}'_t$ is nonempty is eventually periodic, and there exists an eventually quasi-polynomial $\x(t)$ such that $\x(t)\in S'_t$ when it is nonempty. This is the desired $\x(t)$.

To prove that {3b} holds, we induct on $k$. The base case $k=1$ is Property {3} and so is true.
Inductively applying Property 3 to the sentence
\[F(\x,t)\wedge\big(\x\ne \x_1(t)\big)\wedge\cdots\wedge\big(\x\ne \x_{k-1}(t)\big)\]
yields new elements, and this formula has solutions exactly when $\abs{S_t}\ge k$, so the set of such $t$ is eventually periodic, by Property 3.
%

\end{document}